\documentclass[11pt,reqno]{amsart}
\usepackage{amssymb,amscd}

\begin{document}

\let\kappa=\varkappa
\let\eps=\varepsilon
\let\phi=\varphi

\def\Z{\mathbb Z}
\def\R{\mathbb R}
\def\C{\mathbb C}
\def\Q{\mathbb Q}
\def\P{\mathbb P}

\def\OO{\mathcal O}
\def\CP{\C{\mathrm P}}
\def\RP{\R{\mathrm P}}
\def\conj{\overline}
\def\Beta{\mathrm{B}}
\def\O{\Omega}
\def\e{\epsilon}
\def\z{\zeta}

\renewcommand{\Im}{{\mathop{\mathrm{IT}}\nolimits}}
\renewcommand{\Re}{{\mathop{\mathrm{Re}}\nolimits}}
\newcommand{\codim}{{\mathop{\mathrm{codim}}\nolimits}}
\newcommand{\id}{{\mathop{\mathrm{id}}\nolimits}}
\newcommand{\Aut}{{\mathop{\mathrm{Aut}}\nolimits}}

\newtheorem{mainthm}{Theorem}
\renewcommand{\themainthm}{\Alph{mainthm}}
\newtheorem{splitthm}{Theorem}[mainthm]
\newtheorem{thm}{Theorem}[section]
\newtheorem{lem}[thm]{Lemma}
\newtheorem{prop}[thm]{Proposition}
\newtheorem{cor}[thm]{Corollary}

\theoremstyle{definition}
\newtheorem{exm}[thm]{Example}
\newtheorem{rem}[thm]{Remark}

\title[Extension of holomorphic maps]{Extension of holomorphic maps
  between real hypersurfaces of different dimension}
\author{Rasul Shafikov}
\address{Department of Mathematics, the University of Western Ontario,
London, Canada  N6A 5B7}
\email{shafikov@uwo.ca}
\author{Kaushal Verma}
\address{Indian Institute of Science,
         Dept. of Mathematics,
         Bangalore 560012, India}
\email{kverma@math.iisc.ernet.in}

\thanks{R.S. was supported in part by the Natural Sciences
and Engineering Research Council of Canada. K.V. was supported by DST
(India) grant no.: SR/S4/MS-283/05}

\begin{abstract}
It is proved that the germ of a holomorphic map from a real analytic
hypersurface $M$ in $\mathbb C^n$ into a strictly pseudoconvex
compact real algebraic hypersurface $M'$ in $\mathbb C^N$,
$1<n\le N$ extends holomorphically along any path on $M$.
\end{abstract}

\maketitle

%%%%%%%%%%%%%%%% section %%%%%%%%%%%%%%%%%%%%%%%%%%%%%%%%%%%%%%%%%%%%%
\section{Introduction}

In this paper we consider the problem of analytic continuation of a
germ of a holomorphic map sending a real analytic hypersurface into
another such hypersurface in the special case when the target
hypersurface is real algebraic but of higher dimension. Our principal
result is the following.

\begin{thm}\label{main}
Let $M$ be a smooth real-analytic minimal hypersurface in $\C^n$, $M'$
be a compact strictly pseudoconvex real algebraic hypersurface in
$\C^N$, $1<n\le N$. Suppose that $f$ is a germ of a holomorphic map at a
point $p\in M$ and $f(M)\subset M'$. Then $f$ extends as a holomorphic
map along any path on $M$ with the extension sending $M$ to~$M'$.
\end{thm}

In the equidimensional case the problem of analytic continuation of the
germ of a map between real analytic hypersurfaces attracted a lot of
attention (see, for example, \cite{pi1}, \cite{w}, \cite{Vi},
\cite{pi3}, \cite{sh1} and \cite{sh2}). The problem, which originated
in the work of Poincar\'e \cite{po} (generalized later in \cite{t} and
\cite{a}), is related to other fundamental questions in several
complex variables, such as boundary regularity of proper holomorphic
mappings, the theory of CR maps, and classification of domains in
complex spaces (for the latter connection see \cite{Vi}, \cite{pi3},
\cite{sh2}, \cite{ns}).

The situation seems to be more delicate in the case of different
dimensions. The first result of this type is probably due to
Pinchuk \cite{pi2} who proved that a germ of a holomorphic map
from a strictly pseudoconvex real analytic hypersurface $M\subset
\C^n$ into a sphere $S^{2N-1}$, $1<n\le N$, extends holomorphically
along any path on $M$. Just recently Diederich and Sukhov \cite{ds}
proved that the same extension holds if $M$ is weakly pseudoconvex.
Theorem~\ref{main} is a direct generalization of these results
(although our methods are quite different). Further, in the case when
$\dim M =\dim M'$, Theorem~\ref{main} generalizes the result in
\cite{sh1}, where the hypersurface $M$ was assumed to be essentially
finite, a stronger condition than minimality. Other related results
also include various extensions obtained when both $M$ and $M'$ are
algebraic (see e.g. \cite{hu}, \cite{ss}, \cite{ber}, \cite{cms},
\cite{z}, \cite{m} and references therein), which state that under
certain conditions a map between two real algebraic submanifolds (or
even sets) is algebraic, and therefore extends to a dense open subset
of~$\C^n$.

Much like in the equidimensional case, analytic continuation can be
used to prove boundary regularity of holomorphic maps.

\begin{thm}\label{t2} Let $D$ and $D'$ be smoothly bounded domains in
  $\C^n$ and $\C^N$ respectively, $1<n\le N$, $\partial D$ is
  real-analytic, $\partial D'$ is real algebraic, and let $f:D\to D'$
  be a proper holomorphic map. Suppose there exists a point $p\in
  \partial D$ and a neighbourhood $U$ of $p$ such that $f$ extends
  smoothly to $\partial D\cap U$. Then the map $f$ extends
  continuously to $\overline D$, and the extension is holomorphic on a
  dense open subset of $\partial D$. If $D'$ is strictly pseudoconvex,
  then $f$ extends holomorphically to a neighbourhood of~$\overline D$.
\end{thm}

For $n=N$ a similar result is contained in \cite{sh3}.
We note that without the assumption of smooth extension of $f$ {\it
somewhere} on $\partial D$ the conclusion of Theorem~\ref{t2} is false
in general. Indeed, there exist proper holomorphic maps of balls of
different dimension that do not extend even continuously to the
boundary (\cite{lo},\cite{f2}), or that are continuous up to the
boundary but are not of class $C^2$ (\cite{d},\cite{h}). Further,
there exist proper maps $f: \mathbb B^n\to \mathbb B^N$ which are
continuous up to the boundary, and $f(S^{2n-1})=S^{2N-1}$, provided
that $N$ is sufficiently large~(\cite{g}).

On the other hand, if $f$ is known to extend smoothly to {\it all} of
$\partial D$, then $f$ extends holomorphically everywhere on
$\partial D$ according to \cite{cdms} and \cite{mmz}. We use these
results to obtain holomorphic extension of $f$ somewhere on the
boundary of $D$ to start analytic continuation along $\partial
D$. Also without the assumption of algebraicity, Forstneri\v c
\cite{f1} proved that a proper holomorphic map $f: D\to D'$ between
strictly pseudoconvex domains $D\subset \C^n$, $D'\subset \C^N$,
$1<n \le N$, with real analytic boundaries which extends smoothly to
$\partial D$, necessarily extends holomorphically on a dense open
subset of $\partial D$ (this was recently improved in \cite{ps} by
showing that the extension is holomorphic everywhere provided that
$1<n\le N\le 2n$.)

The above stated theorems follow from a more general result asserting
a local extension of the map $f$ as a correspondence. More
precisely, the following holds.

\begin{thm}\label{corr}
Let $M$ (resp. $M'$) be smooth hypersurfaces in $\C^n$ (resp.
$\C^N$), $1<n\le N$, where $M$ is real analytic and minimal,
and $M'$ is compact real algebraic. Suppose $\Sigma\subset M$ is a
connected open set, and $f:\Sigma \to M'$ is a real analytic CR
map. Let $b\in\partial \Sigma$, and $\partial \Sigma \cap M$ be a
smooth generic submanifold. Then there exists a neighbourhood
$U_b\subset \C^n$ of $b$ such that $f$ extends to a holomorphic
correspondence $F: U_b\to \C^N$ with~$F(U_b\cap M)\subset M'$.
\end{thm}

We note that in the context of Theorem~\ref{main} it follows that
$M$ is pseudoconvex, however, in Theorem~\ref{corr} neither $M$
nor $M'$ has to be pseudoconvex. The extension given by
Theorem~\ref{corr} is guaranteed to be single valued if $M'$
satisfies the property that $Q'_{z'}\cap M'=\{z'\}$ near any $z'\in
M'$. In particular this holds if $M'$ is strictly
pseudoconvex (cf. \cite{f1}).

There are no known results when a similar analytic continuation would
hold under the assumption that $M'$ is merely real analytic. The
problem is not well understood even in the equidimensional case, where
it is only known that the germ of a map $f: M\to M'$ extends along any
path on $M$ when both $M$ and $M'$ are strictly pseudoconvex
(\cite{pi1}, \cite{Vi}). The case of different dimensions seems to be
even more difficult.

\medskip

\noindent{\bf Acknowledgment.} The authors would like to thank
Prof.~J.~Merker for numerous remarks concerning the first draft of the
paper, in particular for pointing out the construction of ellipsoids
used in Section~4.1.

%%%%%%%%%%%%%%%% section %%%%%%%%%%%%%%%%%%%%%%%%%%%%%%%%%%%%%%%%%%%%%
\section{Preliminaries}

Let $M$ be a smooth real analytic hypersurface in $\C^n$, $n>1$, $0\in
M$, and $U$ a neighbourhood of the origin. If $U$ is sufficiently
small then  $M\cap U$ can be identified by a real analytic defining
function $\rho(z,\overline z)$, and for every point $w\in U$ we can
associate to $M$ its so-called Segre variety in $U$ defined as
\begin{equation}
Q_w= \{z\in U : \rho(z,\overline w)=0\}.
\end{equation}
Note that Segre varieties depend holomorphically on the variable
$\overline w$. In fact, in a suitable neighbourhood $U={\ 'U}\times
U_n\subset \C^{n-1}\times \C$ we have
\begin{equation}\label{product}
  Q_w=\left \{z=({'z},z_n)\in U: z_n = h({'z},\overline w)\right\},
\end{equation}
where $h$ is a holomorphic function. From the reality condition on the
defining function the following basic properties of Segre varieties
follow:
\begin{equation}\label{segre1}
  z\in Q_w \ \Leftrightarrow \ w\in Q_z,
\end{equation}
\begin{equation}\label{segre2}
  z\in Q_z \ \Leftrightarrow \ z\in M,
\end{equation}
\begin{equation}\label{segre3}
  w\in M \Leftrightarrow \{z\in U: Q_w=Q_z\}\subset M.
\end{equation}
The set $I_w:=\{z\in U: Q_w=Q_z\}$ is itself a complex analytic subset of
$U$. So \eqref{segre3}, in particular, implies that if $M$ does not
contain non-trivial holomorphic curves, then there are only finitely
many points in $U$ that have the same Segre variety (for $U$
sufficiently small). For the proofs of these and other properties of
Segre varieties see e.g. \cite{dw}, \cite{df1} or \cite{ber1}.

A hypersurface $M$ is called {\it minimal} if it does not contain any
germs of complex hypersurfaces. In this case the dimension of the set
$I_w$ can be positive (but less than $n-1$) for all~$w\in M$.

If $f: U \to U'$, $U\subset \C^n$, $U'\subset \C^N$, is a
holomorphic map sending a smooth real analytic hypersurface
$M\subset U$ into another such hypersurface $M'\subset U'$, and $U$ is
as in \eqref{product}, then $f(z)=z'$ implies $f(Q_z)\subset
Q'_{z'}$ for $z$ close to the origin. This invariance property of
Segre varieties will play a fundamental role in our arguments.
We will also denote by $w^s$ the {\it symmetric} point of a point
$w=({'w},w_n)\in U$, which is by definition the unique point defined
by $Q_w\cap \{z\in U: {'z}={'w}\}$.

Suppose now that the hypersurface $M\subset \C^N$ is smooth, compact,
connected, and defined as the zero locus of a real polynomial
$P(z,\overline z)$. Then we may define Segre varieties associated with
$M$ as projective algebraic varieties in $\P^N$. Further, this can be
done for every point in $\P^N$. Indeed, let $M$ be given as a
connected component of the set defined by
\begin{equation}
  \{z\in \C^N : P(z,\overline z)=0\}.
\end{equation}
We can projectivize the polynomial $P$ to define $M$ in $\P^N$ in
homogeneous coordinates
\begin{equation}\label{homcoord}
\hat z=[\hat z_0,\hat z_1,\dots,\hat z_N],\ \
z_k=\frac{\hat z_k}{\hat z_0},\ \ k=1,\dots,N,
\end{equation}
as a connected component of the set defined by
\begin{equation}\label{polys1}
  \{\hat z\in \P^N : \hat P(\hat z,\overline{\hat z})=0 \}.
\end{equation}
We may define now the {\it polar} of $M$ as
\begin{equation}\label{polys2}
  \hat M^c=\{(\hat z,\hat \z)\in \P^N \times \P^N : \hat
  P(\hat z,\hat \z)=0\}.
\end{equation}
Then $\hat M^c$ is a complex algebraic variety in $\P^N\times
\P^N$. Given $\tau\in \P^N$, we set
\begin{equation}\label{polys3}
  \hat Q_\tau = \hat M^c \cap \{(\hat z,\hat\z)\in \P^{N}\times\P^N
  : \hat\z =\overline \tau\}.
\end{equation}
We define the projection of $\hat Q_\tau$ to the first coordinate to be
the Segre variety of~$\tau$.

Recall that for domains $D\subset \C^n$ and $D'\subset \C^N$, a
holomorphic correspondence $F : D \to D'$ is a complex analytic set $A
\subset D \times D'$ of pure dimension $n$ such that the coordinate
projection $\pi: A \to D$ is proper (while $\pi': A\to D'$ need not
be). In this situation, there exists a system of canonical defining
functions
\begin{equation}\label{a}
\Phi_I(z,z')=\sum_{|J|\le m}\Phi_{IJ}(z)z'^J,\ (z,z')\in D\times D',
\ |I|=m,
\end{equation}
where $\Phi_{IJ}(z)$ are holomorphic on $D$, and $A$ is the
set of common zeros of the functions $\Phi_I(z,z')$. For details see,
e.g. \cite{c}. It follows that $\pi$ is in fact surjective and a
finite-to-one branched covering. In particular, there exists a complex
subvariety $S\subset D$ and a number $m$ such that
\begin{equation}
F:=\pi'\circ \pi^{-1}=\{f^1(z),\dots,f^m(z)\},
\end{equation}
where $f^j$ are distinct holomorphic maps in a neighborhood of
$z\in D\setminus S$. The set $S$ is called the {\it branch locus} of
$F$. We say that the correspondence $F$ {\it splits} at $z\in {D}$ if
there is an open subset $U\ni z$ and holomorphic maps $f^j:U\to D'$,
$j=1,2,\dots,m,$ that represent~$F$.

%%%%%%%%%%%%%%%% section %%%%%%%%%%%%%%%%%%%%%%%%%%%%%%%%%%%%%%%%%%%%%
\section{Proof of Theorem \ref{corr}}

In the proof of Theorem \ref{corr} we modify the approach in \cite{sh1}
to our situation. The strategy can be outlined as follows. Without
loss of generality we may assume that $f$ is a holomorphic map
defined in a neighbourhood of $\Sigma$, and $f(\Sigma)\subset
M'$. According to \cite{sh1}, there exists a dense open subset $\omega$
of $Q_b$ with the property that for $a\in \omega$, $Q_a\cap \Sigma \ne
\varnothing$. Furthermore, there exists a non-constant curve
$\gamma\subset \Sigma\cap Q_a$ with the endpoint at $b$. Thus we have
a choice of points $\xi$ and $a$ such that
\begin{equation}\label{points}
a\in Q_b,  \ \ \ \xi\in \gamma\subset\Sigma\cap Q_a.
\end{equation}
The extension of $f$ to the point $b$ can be proved in two
steps. Suppose that $f$ is holomorphic in $U_\xi$, some neighbourhood
of $\xi$. Let $U$ be a neighbourhood of $Q_\xi$. We first show that
the set $A$ defined by
\begin{equation}\label{A}
  A=\left\{(w,w')\in U\times \mathbb C^N : f(Q_w\cap U_\xi)\subset Q'_{w'}
  \right\}
\end{equation}
is complex analytic with the property that $A$ contains $\Gamma_f$,
the graph of $f$, and the projection $\pi: A\to U$ is onto. Further,
$A$ can be extended to an analytic subset of $U\times \P^N$, and we
denote by $\pi': A\to\P^N$ the other coordinate projection.

Secondly, we choose suitable neighbourhoods $U_a$ and $U^*$ of $a$ and
$Q_a$ respectively, and consider the set
\begin{equation}\label{A^*}
  A^*=\left\{(w,w')\in U^*\times \P^N : \pi^{-1}(Q_w\cap
  U_a)\subset \pi'^{-1}(Q'_{w'})
  \right\}.
\end{equation}
We then show that $A^*$ also contains the graph of $f$, and its
projection $\pi^*$ to the first component is onto. In particular,
$\pi^*(A^*)$  contains a neighbourhood of $b$. Note that by
construction the dimension of $A$ may be bigger than $n=\dim
\Gamma_f$. An important fact, however, is that $\dim A^*=n$,
regardless of the dimension of the set $A$. This allows us to show
that $f$ extends locally as a holomorphic correspondence to a
neighbourhood of~$b$.

%%%%%%% subsection 1 %%%%%%%%%%%%%%%%%%%%%%%%%%%%%%%%%%%%%%%%
\subsection{Extension along $Q_\xi$} In this subsection we show that
if $f$ is holomorphic at $\xi\in \Sigma$, then we can extend the graph of
$f$ as an analytic set along $Q_\xi$. It follows from \eqref{segre1}
that there exist neighbourhoods $U_\xi$ of $\xi$ and $U$ of $Q_\xi$
such that for any point $w\in U$, the set $Q_w\cap U_\xi$ is
non-empty. Further, $U_\xi$ and $U$ can be chosen such that $Q_w\cap
U_\xi$ is connected for all $w\in U$. We claim that the
set defined by \eqref{A} is a closed complex analytic subset of
$U\times \C^N$. Indeed, the inclusion $f(Q_w\cap U_\xi)\subset
Q'_{w'}$ can be expressed (cf. \cite{sh1}) as
\begin{equation}\label{system}
  P'\left(f({'z},h({'z},\overline w)),\overline w'\right)=0,
\end{equation}
where $P'(z',\overline z')$ is the defining polynomial of $M'$, and $h$
is the map defined in \eqref{product}. After
conjugation this becomes a system of holomorphic equations in $w$ and
$w'$. The variable $'z$ plays the role of a parameter here, but from
the Noetherian property of the ring of holomorphic functions, we may
extract a finite subsystem which defines $A$ as a complex analytic
set. Further, since the equations in \eqref{system} are polynomials in
$w'$, we may projectivize $A$. This defines an analytic set in
$U\times \P^N$, which we denote again by $A$ for simplicity.

Finally, observe that by the invariance property of Segre
varieties it follows that $A$ contains the points of the form
$(w,f(w))$, $w\in U_\xi$, and therefore $A$ contains the germ at $\xi$
of the graph of $f$. This also shows that $A$ is not empty. We may
consider only the irreducible components of the least dimension which
contain $\Gamma_f$. Thus we may assume that $\dim A\equiv m\ge n$.

%%%%%%%%%% subsection 2 %%%%%%%%%%%%%%%%%%%%%%%%%%%%%%%%%%%%
\subsection{Extension along $Q_a$} Let $\pi: A\to U$ and $\pi': A\to
\mathbb \P^N$ be the natural projections. Since $\P^N$ is compact, and
$A$ is closed in $U\times \mathbb P^N$, the projection $\pi$ is
proper. By the Remmert proper mapping theorem, $\pi(A)$ is a complex
analytic subset of $U$, which simply means that $\pi(A)=U$. For
$\zeta\in A$ let $l_\zeta\pi\subset A$ be the germ of the fibre
$\pi^{-1}(\pi(\zeta))$ at $\zeta$. Then for a generic point $\zeta\in
A$, $\dim l_\zeta\pi = m-n$ which is the smallest possible dimension
of the fibre. By the Cartan-Remmert theorem (see e.g. \cite{l}) the
set
\begin{equation}
S:=\{\zeta\in A : \dim l_{\zeta}\pi> m-n\}
\end{equation}
is complex analytic, and by the Remmert proper mapping theorem
$\pi(S)$ is complex-analytic in $U$. We note that $\dim
\pi(S)<n-1$. This can be seen as follows: if $(m-n)+k$ is the generic
dimension of the fibre over $\pi(S)$, $k>0$, then $\dim S= \dim
\pi(S)+(m-n+k)$. Since $\dim S\le m-1$, $\dim \pi(S)\le n-1-k$, and the
assertion holds.

From the above considerations we conclude that $\pi(S)$ does not
contain $Q_b\cap U$. The sets $U$ and $U_\xi$ defined in Section~3.1
certainly depend on the choice of $\xi$. However, if we vary the point
$\xi$ in $\Sigma$, then the sets defined by \eqref{A} with a different
choice of $\xi$ will coincide on the overlaps and satisfy the
properties stated in Section~3.1. Hence, if $a\in \pi(S)\cap Q_\xi
\cap Q_b$, then we may slightly rearrange points $a\in Q_b$ and
$\xi\in \Sigma\cap Q_a$, and repeat the above constructions (keeping
the same notation), so that $a\notin \pi(S)$.

Let $U_a$ be a neighbourhood of the point $a$ in $U$, so small that
$U_a\cap \pi(S)=\varnothing$. Let $\gamma\subset Q_a\cap \Sigma$ be a
path connecting $\xi$ and $b$. We may choose a neighbourhood $U^*$ of
$\gamma$ (including its endpoints) and $U_a$ in such a way that
$Q_w\cap U_a$ is non-empty and connected for any $w$ in
$U^*$. Consider the set $A^*$ defined in~\eqref{A^*}.

\medskip

\noindent{\it Lemma 3.1. $A^*$ is a complex-analytic subset of
$U^*\times \P^N$. }

\begin{proof}
Let $(w_0,w'_0)\in A^*$ be an arbitrary point. Consider
$\pi^{-1}(Q_{w_0}\cap U_a)$. This is a complex analytic subset of
$A\cap(U_a\times\P^N)$. Since $U_a\cap \pi(S)=\varnothing$, the
fibres of $\pi$ are of constant dimension for points in
$U_a$. Therefore, $\pi^{-1}(Q_{w^0}\cap U_a)$ has constant dimension
$m-1$. It follows that analytic sets $\pi^{-1}(Q_{w}\cap U_a)$ have
the same dimension and vary analytically as $w$ varies near $w_0$.
We denote by $B(X,\epsilon)$ the open $\epsilon$-neighbourhood of a
set~$X$.

Let $q\in \pi^{-1}(Q_{w_0}\cap U_a)$. Then there exists an affine
coordinate patch $U'\subset \P^N$, $q \in U_a\times U'$, with
coordinates
\begin{equation}
(z,\zeta')=(z_1,\dots, z_n,\zeta'_{n+1},\dots,\z'_{n+N})
\in U_a\times U',
\end{equation}
and a choice of a coordinate plane in $U_a\times U'$
passing through~$q$, which is spanned by
\begin{equation}
(z_1,z_2,\dots,z_{n-1},\z'_{k_1},\z'_{k_2}\dots,\z'_{k_{m-n}})
\end{equation}
for some $k_1,k_2,\dots k_{m-n}$, such that for some $\e_q>0$, the set
$\pi^{-1}(Q_{w_0}\cap U_a)\cap B(q, \epsilon_q)$ can be represented as
in~\eqref{a}, i.e. as the zero locus of the functions
\begin{equation}
\Phi_I(z,\z')=\sum_{\substack{0\le j \le m_q\\ |J|\le M_q}}
\Phi_{IjJ}
(z_1,z_2,\dots,z_{n-1},\z'_{k_1},\z'_{k_2}\dots,\z'_{k_{m-n}})
(z_n)^j (\tilde \z')^J, \ |I| \le l_q,
\end{equation}
where $\tilde \z'$ are the remaining $(N-m+n)$ coordinates in $U'$,
$J=(j_1,\dots, j_{N-m+n})$, and $\Phi_{IjJ}$ are holomorphic
functions. Since $\pi^{-1}(Q_w\cap U_a)$ depend anti-holomorphically
on $w$, there exists $\delta_q>0$ and a connected open neighbourhood
$\Omega_q\subset B(q, \epsilon_q)$ of the point $q$, such that for
$|w-w_0|<\delta_q$ a similar representation also holds for
$\pi^{-1}(Q_w\cap U_a)\cap \Omega_q$ with functions
\begin{equation}\label{graph}
\Phi_I(z,\z',\overline w)=\sum_{\substack{0\le j\le m_q\\ |J|\le M_q}}
\Phi_{IjJ} (z_1,\dots,z_{n-1},\z'_{k_1},\dots,\z'_{k_{m-n}},
\overline w) (z_n)^j (\tilde \z')^J, \ |I| \le l_q,
\end{equation}
where the dependence on $\overline w$ is holomorphic.

We claim that there exist $\delta>0$ and a finite collection of points
$q^k\in \pi^{-1}(Q_{w_0}\cap U_a)$, $k=1,2,\dots l$ such that
$\cup_{k=1}^l \Omega_{q^k}$ has a non-empty intersection with every
irreducible component of $\pi^{-1}(Q_{w}\cap U_a)$, provided that
$|w-w_0|<\delta$.

To prove the claim first observe that from compactness of $\P^N$ and
continuity of the fibres of the projection $\pi$, it follows that given
any small $\e>0$ there exists $\delta>0$ such that the distance between
$\pi^{-1}(Q_w\cap U_a)$ and $\pi^{-1}(Q_{w_0}\cap U_a)$ is less than
$\e$ whenever $|w-w_0|<\delta$. The distance in $U_a\times  \P^N$ can
be taken with respect to the product metric of the standard metric in
$\C^n$ and the Fubini-Study metric in~$\P^N$.

Denote by $S^j_{w}$ the irreducible components of $\pi^{-1}(Q_{w}\cap
U_a)$, $j=1,\dots,l_w$, where $w$ is a point in some small
neighbourhood of $w_0$. Choose $\epsilon_1 >0$ and $\delta_1>0$ such
that for $|w-w_0|<\delta_1$, none of the components $S^j_{w}$ is
entirely contained in $B(\partial U_a \times \P^N,\epsilon_1)$. Such
$\epsilon_1$ and $\delta_1$ exist because every $S^j_{w}$ surjectively
projects onto $Q_w\cap U_a$. Then $(U_a\times \P^N )\setminus
B(\partial U_a \times \P^N,\epsilon_1)$ is compact, and therefore, the
open cover of the set
\begin{equation}
\pi^{-1}(Q_{w_0}\cap U_a)\setminus B(\partial U_a \times
\P^N,\epsilon_1)
\end{equation}
by $\Omega_q$, where $q \in \pi^{-1}(Q_{w_0}\cap U_a)$,
admits a finite subcover, say, $\Omega_{q^1},\dots \Omega_{q^l}$. Let
\begin{equation}
\e_2 = \min_{k=1,\dots,l} \left\{\sup \{\alpha>0 :
B(q^k,\alpha)\subset \Omega_{q^k}\} \right\}.
\end{equation}
Then there exists $\delta_2$ such that the
distance between $\pi^{-1}(Q_w\cap U_a)$ and $\pi^{-1}(Q_{w_0}\cap
U_a)$ is less than $\e_2$ whenever $|w-w_0|<\delta_2$. Finally, choose
$\delta=\min\{\delta_1,\delta_2\}$. Then for any $w$ with
$|w-w_0|<\delta$, any component $S^j_{w}$ has a non-empty intersection
with $\cup_k \Omega_{q^k}$. This proves the claim.

We now show that $A^*$ is complex-analytic in a neighbourhood of a
point $(w_0,w'_0)\in A^*$. Choose $q^1,\dots q^l$ as claimed above.
We fix some $q^k$, $k\in \{1,2,\dots,l\}$ and let $\eta=\overline w$,
$\eta'=\overline w'$. Let further
$G=\Omega_{q^k}\times\{|(\eta,\eta')-(\eta_0,\eta'_0)|<\delta\}$
be a small neighbourhood of $(q^k,\overline w_0, \overline w'_0)$ in
$\C_{z}^n\times \C_{\zeta'}^{N} \times \C_{\eta}^n \times \C_{\eta'}^{N}$.
We define
\begin{equation}
X_1 = \{(z,\z',\eta,\eta')\in G: P'(\z',\eta')=0\},
\end{equation}
\begin{equation}
X_2 = \{(z,\z',\eta,\eta')\in G: \Phi^k_I(z,\z',\eta)=0, \ |I| \le l_{q^k}\},
\end{equation}
where $\Phi^k_I(z,\z',\eta)$ are holomorphic functions in as defined 
in \eqref{graph}. Both of these sets are complex analytic in $G$. 
Then the set of points $(w,w')$ for which the inclusion
\begin{equation}\label{inc}
\pi^{-1}(Q_w\cap U_a)\cap \Omega_{q^k}\subset \pi'^{-1}(Q'_{w'})
\end{equation}
holds is conjugate to the set $X^*$ in the $(\eta,\eta')$ space
which is characterized by the property that $(\eta,\eta')\in X^*$
whenever $\pi_2^{-1}(\eta,\eta') \subset \pi_1^{-1}(\eta,\eta')$,
where $\pi_j$ is the coordinate projection from $X_j$ to the
$(\eta,\eta')$-space. The set $X^*$ can be also defined as
\begin{equation}
X^*=\{(\eta,\eta'): \dim \pi_2^{-1}(\eta,\eta') = \dim \pi_{12}^{-1}(\eta,\eta')\},
\end{equation}
where $\pi_{12}: X_1\cap X_2 \to \C^{n+N}_{(\eta,\eta')}$. Further,
$\dim \pi_2^{-1}(\eta,\eta') = m-1$,
for all $(\eta,\eta')$, and so $\dim \pi_{12}^{-1}(\eta,\eta')\le m-1$.
Thus, $X^*=\pi_{12}(\tilde X)$, where
\begin{equation}
\tilde X = \{(z,\z',\eta,\eta')\in X_1\cap X_2: \dim l_{(z,\z',\eta,\eta')}\pi_{12}>m-2\}.
\end{equation}
By the Cartan-Remmert theorem $\tilde X$ is a complex analytic subset of $G$. Denote
by $\tilde \pi$ the projection from $\tilde X$ to the space of variables
$(z_1,\dots,z_{n-1},\zeta_1,\dots,\z_{k_{m-n}},\eta,\eta')$. By construction of
functions in \eqref{graph} the map $\tilde \pi$ is proper. Hence, by the Remmert proper mapping
theorem, $\tilde \pi (\tilde X)$ is complex analytic. Finally, consider the
projection $\pi_{(\eta,\eta')}: \tilde\pi(\tilde X)\to (\eta,\eta')$. From
the construction of the set $\tilde\pi(\tilde X)$,
$\dim \pi_{(\eta,\eta')}^{-1}(\eta,\eta')=m-1$, for $(\eta,\eta')\in X^*$.
But in fact, $\dim \pi_{(\eta,\eta')}^{-1}(\pi_{(\eta,\eta')}(x))=m-1$, for any
$x\in \tilde\pi(\tilde X)$. Thus we may identify $X^*$ with
$\tilde\pi(\tilde X)\cap\{(z_1,\dots,z_{n-1},\zeta_1,\dots,\z_{k_{m-n}})={\rm const}\}$.
This proves that the set $X^*$ is complex analytic.
After conjugation, we may assume that the set defining the inclusion
in \eqref{inc} is also complex analytic.

If an open set of the irreducible component $S^j_{w}$ is
contained in $\pi'^{-1}(Q'_{w'})$ for some $w'$, then by the
uniqueness theorem, the whole component $S^j_{w}$ must be contained in
$\pi'^{-1}(Q'_{w'})$. Therefore, since $\cup_{k=1}^l \Omega_{q^k}$ has
a non-empty intersection with every $S^j_w$, the system of equations
defining the inclusion \eqref{inc}, combined for $k=1,\dots, l$,
completely determines the inclusion in \eqref{A^*}, and therefore it
defines $A^*$ as a complex-analytic set near~$(w_0,w'_0)$.

So far we have showed that $A^*$ is a {\it local} complex analytic
set, i.e defined by a system of holomorphic equations in a
neighbourhood of any of its points. To prove that $A^*$ is a
complex-analytic {\it subset} of $U^*\times \P^N$ it is enough
now to show that $A^*$ is closed in $U^*\times \P^N$. Suppose
$(w^j,w'^j)\to (w^0,w'^0)$, as $j\to\infty$, for some sequence
$(w^j,w'^j)\in A^*$, and suppose that $(w^0,w'^0)\in U_a\times
\P^N$. This means that $\pi^{-1}(Q_{w^j}\cap
U_a)\subset\pi'^{-1}(Q'_{w'^j})$. Since $Q_{w^j}\to Q_{w^0}$ and
$Q'_{w'^j}\to Q'_{w'^0}$, by analyticity also
$\pi^{-1}(Q_{w^0}\cap U_a)\subset \pi'^{-1}(Q'_{w'^0})$, and
therefore $(w^0,w'^0)\in A^*$. This completes the proof of Lemma
3.1.
\end{proof}

\noindent{\it Lemma 3.2. The set $A^*$ contains the germ of the graph
  of $f$ at $(\xi, f(\xi))$. Further, $$A^*\cap \left((U_\xi \cap
  U\cap U^*)\times \mathbb P^N \right) \subset A.$$}

\begin{proof} Suppose $z\in (U_\xi \cap U\cap U^*)$. We need
to show that
\begin{equation}\label{z}
\pi^{-1}(Q_z\cap U_a)\subset \pi'^{-1}\left(Q'_{f(z)}\right).
\end{equation}
Let $w\in Q_z\cap U_a$ be an arbitrary point, and let $(w,w')\in
A$. Then $f(Q_w\cap U_\xi)\subset Q'_{w'}$. In particular, since $z\in
Q_w\cap U_\xi$, we have $f(z)\in Q'_{w'}$. But this implies $w'\in
Q'_{f(z)}$. In other words, $(w,w')\in \{w\}\times Q'_{f(z)}$. Since
$w'$ was an arbitrary point in $A$ over $w$, we conclude that
$\pi^{-1}(w)\subset \pi'^{-1}\left(Q'_{f(z)}\right)$. Consequently,
\eqref{z} follows, and $(z,f(z))\in A^*$.

As for the second assertion, we observe that for $(w,w')\in A^*$,
where~$w$ is sufficiently close to~$\xi$, the inclusion
$\pi^{-1}(Q_w\cap U_a)\subset \pi'^{-1}(Q'_{w'})$ is equivalent to
$\pi^{-1}(Q_w\cap U_\xi)\subset \pi'^{-1}(Q'_{w'})$, because $Q_w\cap
U$ is connected. From Section~3.1 the set $A$ contains the germ of the
graph of $f$ near $\xi$, and therefore, the inclusion $\Gamma_f\subset
A^*$ in particular implies $f(Q_w\cap U_\xi)\subset Q'_{w'}$, which by
definition means $(w,w')\in A$.
\end{proof}

Lemma~3.2 shows that $A^*$ is non-empty. Also note
that since $\P^N$ is compact, the projection $\pi^*:A^*\to U^*$ is
proper, and therefore, $\pi^*(A^*)=U^*$. Define $\pi'^*: A^*\to
\P^N$.

%%%%%%%%%% subsection 3 %%%%%%%%%%%%%%%%%%%%%%%%%%%%%%%%%%%%
\subsection{Extension as a correspondence.}

Let now $\O$ be a small connected neighbourhood of the path
$\gamma\subset Q_a\cap M$, which connects~$\xi$ and~$b$, such
that for any $w\in \O$, the symmetric point~$w^s$ belongs to $U^*$,
and let $Q^s_w$ denote the connected component of $Q_w\cap U^*$ which
contains~$w^s$. Denote further by $S^*$ the set of points $z\in
U^*$ for which ${\pi^*}^{-1}(z)\subset A^*$ does not have the generic
dimension. The same argument as at the beginning of Section~3.2 shows
that $S^*$ is a complex analytic set of dimension at most $n-2$, and
so $\O\setminus S^*$ is connected. To prove the extension of $f$ to
the point $b$ we will need the following result.

\medskip

\noindent{\it Lemma~3.3. For any point $w\in \O\setminus S^*$,
\begin{equation}\label{i}
{\pi^*}^{-1}(Q^s_w)\subset {\pi'^*}^{-1}(Q'_{w'}), \ \ \forall\
w'\in\pi'^*\circ{\pi^*}^{-1}(w).
\end{equation}}

\begin{proof}
Denote by $Z$ the set of points in $\O\setminus S^*$ for which
\eqref{i} holds. We show that $Z=\O\setminus S^*$. For the proof we
shrink $U_\xi$ so that $U_\xi \subset \Omega$.

Let $w\in U_\xi \setminus S^*$ be some point, and $(w,w')\in
A^*$. Note that $z^s=z$ for any $z\in M$, and therefore, for $w$
sufficiently close to $\xi$, the set $Q_w\cap U_\xi$ coincides with
$Q^s_w\cap U_\xi$. Let $z\in Q_w\cap U_\xi$ be arbitrary. Then
$(z,z')\in A^*$ means  $\pi^{-1}(Q_z\cap U_a)\subset
\pi'^{-1}(Q'_{z'})$. For $z$ and $w$ sufficiently close to $\xi$,
$Q_z$ is connected in $U$, and therefore, $\pi^{-1}(Q_z\cap
U_\xi)\subset \pi'^{-1}(Q'_{z'})$. The last inclusion in particular
means that $\pi^{-1}(w) \subset\pi'^{-1}(Q'_{z'})$. Thus for any
$w'\in \pi'\circ\pi^{-1}(w)$, $w'\in Q'_{z'}$, or $z'\in Q'_{w'}$. By
Lemma~3.2, $A^*$ is contained in $A$ near $\xi$, and it follows that
for any $w'\in{\pi^*}'\circ{\pi^*}^{-1}(w)$, $z'\in Q'_{w'}$. From
that \eqref{i} follows, and we proved that the set $Z$ contains a
small neighbourhood of~$\xi$.

Let $Z^{\circ}$ be the largest connected open set which contains $\xi$
and is contained in $Z$. From the above considerations, $Z^\circ\ne
\varnothing$. We show that if $w \in (\overline{Z^\circ} \setminus
Z^\circ)\cap (\Omega\setminus S^*)$, then $w\in Z^\circ$. Let
$(w,w')\in A^*$ for some $w'$. Since $\dim S^* < \dim Q^s_w=n-1$, we
may find a point $\alpha\in (Q^s_w\setminus S^*)$, and by repeating
the argument of  Lemma~3.1 we may construct a complex analytic set
\begin{equation}\label{xx'}
A_w=\{(x,x')\in U_w\times \P^N:
{\pi^*}^{-1}(Q^s_x\cap U_\alpha) \subset {{\pi^*}'}^{-1}(Q'_{x'})\},
\end{equation}
where $U_w$ and $U_\alpha$ are suitably chosen neighbourhoods of $w$
and $\alpha$ respectively. For every point $x\in U_w\cap Z^\circ$, and
every $x'$ such that $(x,x')\in A^*$, the inclusion in \eqref{xx'}
holds. This implies
\begin{equation}
A^*\cap \left((Z^\circ \cap U_w)\times\P^N \right)\subset A_w,
\end{equation}
and in particular, $A_w$ is non-empty. By the uniqueness theorem,
it follows that $A^*\cap (U_w \times\P^N)\subset A_w$, and therefore,
the projection from $A_w$ to the first component is onto. Thus, for
any $x\in U_w$, the set $Q_x\cap U_\alpha$ (and therefore $Q^s_x$)
will be ``mapped'' by $A^*$ into Segre variety of a point $x'$,
whenever $(x,x')\in A^*$. Hence, $U_w\subset Z^\circ$.

Since $\Omega\setminus S^*$ is connected, it follows now that
$Z=\Omega\setminus S^*$.
\end{proof}

We now consider only an irreducible component of $A^*$ which has the
smallest dimension, and such that it contains the germ of the graph of
$f$ at $\xi$. Denote for simplicity this component again by
$A^*$. Note that Lemma~3.3 still holds for the new $A^*$.

\medskip

\noindent{\it Lemma 3.4. $\dim A^*=n$.}

\begin{proof}
Since $\pi^*: A^*\to U^*$ is onto, for any $z\in M\cap U^*$ the set
${\pi^*}^{-1}(z)$ is non-empty. We show that for a given $z_0\in
\Sigma\setminus S^*$, the set ${\pi^*}^{-1}(z_0)$ is
discrete near $(z_0,f(z_0))\in A^*$. Indeed, by Lemma~3.3,
$(z,z')\in A^*\setminus {\pi^*}^{-1}(S^*)$ implies
${\pi^*}^{-1}(Q^s_{z})\subset {\pi'^*}^{-1}(Q'_{z'})$. In particular
this means that $\pi'^*({\pi^*}^{-1}(z))\subset Q'_{z'}$, which
implies that $z'\in Q'_{z'}$. Then from~\eqref{segre2} it follows that
for any $z\in M$ close to $z_0$, and any $z'$ close to $f(z_0)$, the
inclusion $(z,z')\in A^*$ implies $z'\in M'$. Since
$\pi'^*({\pi^*}^{-1}(z))$ is a locally countable union of complex
analytic sets, and $M'$ contains no non-trivial germs of complex
analytic varieties by \cite{df}, it follows that ${\pi^*}^{-1}(z)$ is
discrete near $(z_0,f(z_0))$. This means that $\dim A^*=n$ near $(z_0,
f(z_0))$. But then the lemma follows, since $\dim A^*$ is
constant.
\end{proof}

To finish the proof of the theorem, we consider two cases. First,
suppose that $b\not\in S^*$. Since $M'$ is compact, the cluster set of
$f|_\gamma (b)$ is well-defined. Let $b'$ be a point in the cluster
set of the point $b$. It is enough to show now that there exist
neighbourhoods $U_b\ni b$ and $U'_{b'}\ni b'$ such that $A^*\cap
U_b\times U'_{b'}$ is a holomorphic correspondence. By construction,
$(b,b')\in A^*$, and from the proof of Lemma~3.4 we conclude that
${\pi^*}^{-1}(b)$ is discrete near $(b,b')$. Therefore we may choose
$U_b$ and $U'_{b'}$ in such a way that $A^*\cap (U_b\times \partial
U'_{b'})=\varnothing$. It follows then that $\pi^*|_{A^*\cap
(U_b\times U'_{b'})}$ is a proper map, and therefore,
\begin{equation}\label{F}
F:={\pi'^*}|_{A^*\cap (U_b\times U'_{b'})}\circ{\pi^*}^{-1}|_{U_b}
\end{equation}
is the desired extension of $f$ as a holomorphic correspondence.

Secondly, suppose $b\in S^*$. Consider a sequence of points $w^j\in
(\Sigma\cap \Omega)\setminus S^*$ such that $w^j \to b$ and $\lim
f(w^j)=b'$ for some $b'\in M'$. Then
\begin{equation*}
{\pi^*}^{-1}(Q^s_{w^j})\subset {\pi'^*}^{-1}(Q'_{{f(w^j)}}).
\end{equation*}
It follows that
\begin{equation}\label{bb'}
{\pi^*}^{-1}(Q^s_{b})\subset Q'_{{b'}}.
\end{equation}
Indeed, it is enough to prove this inclusion in a neighbourhood of any
point in $Q^s_b$. Since $\dim S^* < \dim Q_b$, we may choose this
point to be outside $S^*$. The inclusion then follows by analyticity
of the fibres of $\pi^*: A^*\to U^*$.

As in the proof of Lemma 3.4, it follows from \eqref{bb'} that
${\pi^*}^{-1}(b)$ is discrete near $(b,b')$, and the same argument as
above shows that $f$ extends to a neighbourhood of $b$ as a
holomorphic correspondence.

Finally, if $F$ is the extension of $f$ as a correspondence, then
$F(M)\subset M'$. The reason again is that if $z\in M$ and $z'\in
F(z)$, then $F(Q_z)\subset Q'_{z'}$ by \eqref{i} and \eqref{bb'},
which implies $z'\in M'$, and by~\eqref{segre2}, $z'\in M'$.

This completes the proof of Theorem~\ref{corr}.

%%%%%%%%%%%%% section %%%%%%%%%%%%%%%%%%%%%%%%%%%%%%%%%%%%%%%%
\section{Proof of other results.}

\subsection{Proof of Theorem~\ref{main}}

We first show that the map $f$ can be extended holomorphically along
any smooth CR-curve $\gamma$ on $M$, i.e. for which the tangent
vector to $\gamma$ at any point is contained in the complex tangent to
$M$. For this we use the construction of a family of ellipsoids
first used by Merker and Porten \cite{mp}. We refer to their paper for
the details of this construction. Let $q$ be the first point
on $\gamma$ to which $f$ does not extend holomorphically. Near $q$
there exists a smooth CR vector field $L$ such that $\gamma$ is
contained in an integral curve of $L$. By integrating $L$ we obtain a
smooth coordinate system $(t,s)\in \R\times \R^{2n-2}$ on $M$ such
that for any fixed $s_0$ the segments $(t,s_0)$ are contained in the
trajectories of $L$. We may further choose a point $p\in \gamma$
sufficiently close to $q$, so that $f$ is holomorphic near $p$. After
a translation, assume that $p=(0,0)$. For $\e>0$ define the
family of ellipsoids on $M$ by
\begin{equation}\label{e}
E_\tau = \{(t,s) : |t|^2/\tau+|s|^2 < \e\},
\end{equation}
where $\e>0$ is so small that for some $\tau_0>0$ the ellipsoid
$E_{\tau_0}$ is compactly contained in the portion of $M$ where $f$ is
holomorphic. Then $\partial E_\tau$ is generic at every point except
the set
\begin{equation*}
\Lambda=\{(0,s): |s|^2=\e\}.
\end{equation*}
Let further $\tau_1>0$ be such that $q\in \partial E_{\tau_1}$.

To prove that $f$ extends holomorphically to a neighbourhood of $q$ we
argue by contradiction. For that we assume that $\tau^*$ is the
smallest positive number such that $f$ does not extend holomorphically
to {\it some} point on $\partial E_{\tau^*}$, and assume that
$\tau^*<\tau_1$. By construction, $\tau^*>\tau_0$. Also by
construction, near any point $b\in \partial E_{\tau^*}$ to which $f$
does not extend holomorphically, the set $\partial E_{\tau^*}$ is a
smooth generic submanifold of $M$, since the non-generic points of
$\partial E_{\tau^*}$ are contained in $\Lambda$, where $f$ is already
known to be holomorphic. Then by Theorem~\ref{corr} the map $f$
extends as a correspondence $F$ to a neighbourhood of~$b$.

We now show that $F$ is single valued. Suppose  $w'\in F(w)$ for $w\in
M$, then by the invariance of Segre varieties $F(Q_w)\subset Q'_{w'}$,
and in particular, $w'\in Q'_{w'}$. But since $M'$ is strictly
pseudoconvex, in a sufficiently small neighbourhood of $w'\in M'$
there exists only one point on $M'$ whose Segre variety contains $w'$,
namely $Q'_{w'}$ itself. Thus the correspondence $F$ splits into
several holomorphic maps, one of which by analyticity extends the
map~$f$.

This shows that $\tau^*$ cannot be smaller than $\tau_1$, which proves
that the map $f$ extends holomorphically to $q$, and therefore along
any CR-curve on~$M$.

Finally, observe that minimality of $M$ implies that CR-orbit of any
point on $M$ coincides with $M$. Therefore, using analytic
continuation along CR-curves we obtain continuation of $f$ to a
neighbourhood of any point on~$M$.

\subsection{Proof of Theorem~\ref{t2}}
By \cite{cdms} and \cite{mmz} it follows that smooth extension of
$f$ implies holomorphic extension to a neighbourhood of $p$. We may
apply now Theorem~\ref{corr} to obtain extension as a holomorphic
correspondence, which splits into holomorphic maps on a dense open
subset of $\partial D$. Therefore we may use, for example, one of the
branches of the extension to apply Theorem~\ref{corr} again. Clearly
at each step the extension coincides with $f$ in~$D$.

The extension of $f$ as a correspondence in particular proves
continuous extension of $f$ to $\partial D$. Indeed, if $F$
extends $f$ as a correspondence in a neighbourhood of $q \in
\partial D$, then the cluster set of $q$ with respect to $f$ must be
contained in the set $F(q)$ which is finite. Since the cluster
set is connected it must reduce to a single point thereby showing
that $f$ is continuous at $q$. Holomorphic extension on a dense
open subset of $\partial D$ now simply follows from the splitting
property of correspondences.

The second statement of the theorem follows immediately from
Theorem~\ref{main}.

\bigskip

%%%%%%%%%%%%%%%%%%%%% bibliography %%%%%%%%%%%%%%%%%%%%%%%%%%%%%%%%%%%%%%

\end{document}